# A note on an expansion formula with application to nonlinear DAE's

Matthias Stiefenhofer

ABSTRACT. In [DL] systems of differential polynomials are investigated with respect to properties of Artin approximation type. The key tool in [DL] is an extended version of a formula by Hurwitz [Hu] expressing high order derivatives of an expansion by lower ones. The formula is further refined in [VFZ] to deliver sufficient conditions concerning the existence of power series solutions of scalar algebraic differential equations of order $n$.

In the paper at hand, the main results from [VFZ] are first reproduced and further extended to systems of nonlinear differential algebraic equations. In addition, a simple extension of Tougeron's implicit function theorem is given in a specific constellation.

The results follow from [S1], [S2] where Artin approximation is treated within a Banach space setting, thereby constructing an expansion formula that expresses accurately the required dependency of low and high order derivatives within the system of undetermined coefficients.

**Keywords:** Differential polynomials, Expansion formula, Separant matrix, Tougeron implicit function theorem, Nonlinear differential algebraic equations

*Contents*

## *1. An expansion formula*

Power series mappings of the form

$$G[z] = \sum_{i=0}^{\infty} \frac{1}{i!} G_0^i \cdot (z - z_0)^i \tag{1.1}$$

with symmetric multilinear operators $G_0^i \in L[B^i, \bar{B}]$ are investigated in [S1], [S2], [S3] with respect to power series solutions $G[z(\varepsilon)] = 0$ of the form

$$z(\varepsilon) = \sum_{i=0}^{\infty} \frac{1}{i!} z_i \cdot \varepsilon^i \, , \;\; \varepsilon \in \mathbb{K} = \mathbb{R}\, , \mathbb{C} \tag{1.2}$$

with focus layed upon the existence of convergent power series solutions $z(\varepsilon)$ ascertained by the implicit function theorem. $B$ and $\bar{B}$ denote Banach spaces over $\mathbb{K} = \mathbb{R}$ or $\mathbb{K} = \mathbb{C}$.

In some more detail, the ansatz (1.2) is plugged into (1.1), implying a mapping between power series according to

$$z(\varepsilon) = \sum_{i=0}^{\infty} \frac{1}{i!} z_i \cdot \varepsilon^i \quad \longrightarrow \quad G[\, z(\varepsilon)\,] = \sum_{i=0}^{\infty} \frac{1}{i!} T^i(z_i, \dots, z_0) \cdot \varepsilon^i, \tag{1.3}$$

where the key lemma concerning the structure of the system of undetermined coefficients $T^i(z_i, \dots, z_0) = 0, i \in \mathbb{N}$, is given by the following result.

**Lemma 1:** For $k \geq 0,\ l \geq 0$, the coefficient $T^{2k+1+l}(\cdot)$ of $G[\, z(\varepsilon)\,]$ can be expressed according to

$$T^{2k+1+l}(\, z_{2k+1+l}, \dots, z_{k+1+l},\ \ z_{k+l}, \dots, z_{k+1},\ \ z_k, \dots, z_0\,)$$

$$= \left[\, T^0_{z_0}(z_0)\ \ T^2_{z_1}(z_1, z_0)\ \cdots\ T^{2k}_{z_k}(z_k, \dots, z_0)\,\right] \cdot \boxed{C^{2k+l}} \cdot \begin{pmatrix} z_{2k+1+l} \\ \vdots \\ z_{k+1+l} \end{pmatrix} + R_{2k+1+l}(z_{k+l}, \dots, z_0) \tag{1.4}$$

with

$$\boxed{C^{2k+l}} = Diag\left[\, \gamma_0^{2k+l}\ \cdots\ \gamma_k^{2k+l}\,\right] \ \ and\ \ \gamma_t^{2k+l} := \binom{2k+1+l}{t} \cdot \binom{2t}{t}^{-1},\ \ t = 0, \dots, k. \tag{1.5}$$

Lemma 1 summarizes formulas (4.3)-(4.8) in [S2]. Note that Lemma 1 is valid with respect to general Banach spaces $B$ and $\bar{B}$ over $\mathbb{K}$. Generalizations of Lemma 1 are possible.

The coefficients $T^i(z_i, \dots, z_0)$ of the image series in (1.3) are composed of multilinear mappings $G_0^i \in L[B^i, \bar{B}]$ applied to the coefficients $z_0, z_1, \dots$ of the input series . In a finite dimensional setting, the coefficients $T^i(z_i, \dots, z_0)$ are given by polynomials with respect to $z_0, z_1, \dots$ .

In [S1], [S2], formulas (1.4), (1.5) are used to investigate recursively the system of undetermined coefficients

$$\begin{aligned}
T^1(z_1, z_0) \quad &= \overbrace{G'[z_0]}^{=G_0^1} \cdot z_1 = 0 \\
T^2(z_2, z_1, z_0) \quad &= G'[z_0] \cdot z_2 \ + \overbrace{G''[z_0]}^{=G_0^2} \cdot z_1^2 = 0 \\
T^3(z_3, z_2, z_1, z_0) \ &= G'[z_0] \cdot z_3 \ + 3G''[z_0] \cdot z_1 z_2 + \overbrace{G'''[z_0]}^{=G_0^3} \cdot z_1^3 = 0 \\
\vdots \qquad\qquad & \quad \vdots \qquad\qquad\qquad \vdots
\end{aligned} \tag{1.6}$$

obtained by Taylor expansion of $G[z_0 + \varepsilon z_1 + \cdots] = 0,\ G \in C^\infty(B, \bar{B})$, thereby constructing a linear sum operator with increasing range, finally allowing the application of the implicit function theorem if surjectivity of the sum uperator can be reached. The surjectivity condition is closely related to a direct sum condition of order $k$ derived in [ELG], [LG] in the context of local and global bifurcation theory as well as sign change of Brouwer's degree.

The general representation formulas of $T^{2k+1+l}(\cdot)$ by derivatives of $G[z]$ at $z = z_0$, as in (1.6), can be found in detail in [S1] and [S2]. In addition, some of the results are applied to singular perturbation problems of nonhyperbolic points and periodic orbits of ordinary differential equations [BS], [S4], [S5]. In [S3], the results are used to investigate the family of solution curves $z(\varepsilon, p)$ of $G[z] = 0$ with respect to stability, regularity and uniqueness. Finally, some relations to arc space $X_\infty$ and Newton polygons as well as the Milnor number are established.

If $k \geq 0$ is fixed, then the coefficients $T^i(\cdot)$ in (1.3) with $i > 2k$, i.e. $T^{2k+1+l}(\cdot)$ with $l \geq 0$, are characterized by the following properties according to (1.4), (1.5).

- A fixed leading linear operator independent of $l$ (red). In some more detail, the operator is given by a sum operator defined by partial derivatives $T^{2i}_{z_i}(z_i, \dots, z_0),\ i = 0, \dots, k$ of $k+1$ low order coefficients $T^0(\cdot), T^2(\cdot), \dots, T^{2k}(\cdot)$ with derivatives taken with respect to low $z$-coefficients $z_0, \dots, z_k$ respectively. In this sense, high order derivatives $T^{2k+1+l}(\cdot)$ within the power series expansion of $G[\, z(\varepsilon)\, ]$ in (1.3) are determined by lower ones $T^0(\cdot), T^2(\cdot), \dots, T^{2k}(\cdot)$ using some sort of recursion formula. Hence, Lemma 1 is a result in the spirit of [Hu], [DL], Lemma 2.2 and [VFZ], Theorem 3.5, 3.6.
- The sum operator by itself only depends on $k+1$ lowest $z$-coefficients $(z_k, \dots, z_0)$ of ansatz (1.2). On contrary, the linearity is given with respect to $k+1$ highest $z$-coefficients $(z_{2k+1+l}, \dots, z_{k+1+l})$, i.e. a clear separation between low order (red) and high order (green) z-coefficients is obtained.
- The map $R_{2k+1+l}(\cdot)$ is adding up the remaining summands within $T^{2k+1+l}(\cdot)$ with decisive property not to depend on $z$-coefficients of order higher than $(k+l)$, i.e. if the coefficients $(z_{k+l}, \dots, z_0)$ are fixed, then the equation

$$T^{2k+1+l}(z_{2k+1+l}, \dots, z_{k+1+l},\ z_{k+l}, \dots, z_0) = 0$$

  of the system of undetermined coefficients is linear with respect to $k+1$ highest $z$-coefficients $(z_{2k+1+l}, \dots, z_{k+1+l})$. This property is crucial when building up the linear sum operator in [S1], [S2].
- Finally, the diagonal matrix $C^{2k+l}$ is defined in a rather simple way by binomial coefficients depending both from $k$ and $l$.

## *2. Nonlinear scalar DAE's*

As a preliminary, note that if the Banach space $B$ is given by a product space of several $\mathbb{K}$-Banach spaces, e.g. $z = (x, y_0, \dots, y_n) \in X \times Y_0 \times \cdots \times Y_n$, then (1.2) reads

$$z(\varepsilon) = \sum_{i=0}^{\infty} \frac{1}{i!} \underbrace{\begin{pmatrix} x_i \\ y_{0,i} \\ \vdots \\ y_{n,i} \end{pmatrix}}_{=z_i} \cdot \varepsilon^i \tag{2.1}$$

yielding

$$T^{2i}_{z_i}(z_i, \dots, z_0) = \left[\, T^{2i}_{x_i}\ T^{2i}_{y_{0,i}}\ \cdots\ T^{2i}_{y_{n,i}}\, \right](z_i, \dots, z_0) \in L[\, X \times Y_0 \times \cdots \times Y_n, \bar{B}\, ]$$

and implying for $k \geq 0,\ l \geq 0$ by Lemma 1 the representation

$$T^{2k+1+l}(z_{2k+1+l},\dots,z_{k+1+l},\; z_{k+l},\dots,z_{k+1},\; z_k,\dots,z_0) =$$

$$\left[\left[T^0_{x_0} \cdots T^0_{y_{n,0}}\right]\left[T^2_{x_1} \cdots T^2_{y_{n,1}}\right]\cdots\left[T^{2k}_{x_k} \cdots T^{2k}_{y_{n,k}}\right]\right]\cdot\boxed{C^{2k+l}}\cdot\begin{pmatrix}\boxed{\begin{matrix}x_{2k+1+l}\\ y_{0,2k+1+l}\\ \vdots\\ y_{n,2k+1+l}\end{matrix}}\\ \vdots\\ \boxed{\begin{matrix}x_{k+1+l}\\ y_{0,k+1+l}\\ \vdots\\ y_{n,k+1+l}\end{matrix}}\end{pmatrix} + R_{2k+1+l}(z_{k+l},\dots,z_0)$$

$$= \left[T^0_{x_0} \cdots T^{2k}_{x_k}\right]\cdot\boxed{C^{2k+l}}\cdot\begin{pmatrix}x_{2k+1+l}\\ \vdots\\ x_{k+1+l}\end{pmatrix} + \left[T^0_{y_{0,0}} \cdots T^{2k}_{y_{0,k}}\right]\cdot\boxed{C^{2k+l}}\cdot\begin{pmatrix}y_{0,2k+1+l}\\ \vdots\\ y_{0,k+1+l}\end{pmatrix} \tag{2.2}$$

$$+ \cdots + \left[T^0_{y_{n,0}} \cdots T^{2k}_{y_{n,k}}\right]\cdot\boxed{C^{2k+l}}\cdot\begin{pmatrix}y_{n,2k+1+l}\\ \vdots\\ y_{n,k+1+l}\end{pmatrix} + R_{2k+1+l}(z_{k+l},\dots,z_0).$$

In the next step, our aim is to use (2.2) in the context of scalar differential algebraic equations of the form

$$G[\,x,y,y',\dots,y^{(n)}\,] = 0, \quad n \ge 0 \tag{2.3}$$

with $B = X \times Y_0 \times \cdots \times Y_n = \mathbb{K}^{2+n}$ and $G \in C^\infty(\mathbb{K}^{2+n},\mathbb{K})$. We are looking for power series solutions (2.1), thereby restricting the coefficients $z_i$ according to

$$z(\varepsilon) = \sum_{i=0}^{\infty}\frac{1}{i!}\,\varepsilon^i\cdot\underbrace{\begin{pmatrix}x_i\\ y_{0,i}\\ \vdots\\ y_{n,i}\end{pmatrix}}_{=z_i} = \underbrace{\begin{pmatrix}0\\ c_0\\ \vdots\\ c_n\end{pmatrix}}_{=z_0} + \varepsilon\cdot\underbrace{\begin{pmatrix}1\\ c_1\\ \vdots\\ c_{n+1}\end{pmatrix}}_{=z_1} + \frac{1}{2}\varepsilon^2\cdot\underbrace{\begin{pmatrix}0\\ c_2\\ \vdots\\ c_{n+2}\end{pmatrix}}_{=z_2} + \cdots, \tag{2.4}$$

i.e. the differential dependency within $z = (x,y,y',\dots,y^{(n)})$ is respected and $x = \varepsilon$ represents the independent variable of (2.3) by $x_1 = 1, x_i = 0,\ i \neq 1$. Further, the existence of an appropriate initial value of (2.3) satisfying

$$G[\,z_0\,] = G[\,0,c_0,\dots,c_n\,] = 0$$

is assumed, yielding by Lemma 1 and (2.2)

$$T^{2k+1+l}(c_{2k+1+l+n},\dots,c_0)$$

$$= \underbrace{\left[T^0_{x_0} \cdots T^{2k}_{x_k}\right]\cdot\boxed{C^{2k+l}}\cdot\begin{pmatrix}x_{2k+1+l}\\ x_{2k+l}\\ \vdots\\ x_{k+1+l}\end{pmatrix}}_{=\begin{cases}T^0_{x_0} & if\ k=l=0\\ 0 & otherwise\end{cases}} + \left[T^0_{y_{0,0}} \cdots T^{2k}_{y_{0,k}}\right]\cdot\boxed{C^{2k+l}}\cdot\begin{pmatrix}c_{2k+1+l}\\ c_{2k+l}\\ \vdots\\ c_{k+1+l}\end{pmatrix} \tag{2.5}$$

$$+ \cdots + \left[T^0_{y_{n,0}} \cdots T^{2k}_{y_{n,k}}\right]\cdot\boxed{C^{2k+l}}\cdot\begin{pmatrix}c_{2k+1+l+n}\\ c_{2k+l+n}\\ \vdots\\ c_{k+1+l+n}\end{pmatrix} + R_{2k+1+l}(c_{k+l+n},\dots,c_0).$$

Note that the operators in square brackets (red) merely depend on $(z_k, \dots, z_0)$, i.e. by $(c_{k+n}, \dots, c_0)$ under consideration of (2.4). Additionally, the first summand equals zero except in case of $k = l = 0$.

Now, the coefficients $c_{2k+l+n}, \dots, c_{k+1+l+n}$ of the green last vector also occur in previous green vectors, where combining same coefficients implies by elemetary calculation

$$
\begin{aligned}
&T^{2k+1+l}(c_{2k+1+l+n}, \dots, c_0) \\
&= \underbrace{\begin{bmatrix} \gamma_0^{2k+l} & \cdots & \gamma_k^{2k+l} \end{bmatrix}}_{=:\Gamma_k(l)} \cdot \underbrace{\begin{pmatrix} T^0_{y_{n,0}} & T^0_{y_{n-1,0}} & T^0_{y_{n-2,0}} & \cdots & T^0_{y_{n-k,0}} \\ & T^2_{y_{n,1}} & T^2_{y_{n-1,1}} & \cdots & T^2_{y_{n-k+1,1}} \\ & & T^4_{y_{n,2}} & \cdots & T^4_{y_{n-k+2,2}} \\ & & & \ddots & \vdots \\ & & & & T^{2k}_{y_{n,k}} \end{pmatrix}}_{=:S_k(c_{k+n},\dots,c_0)} \cdot \begin{pmatrix} c_{2k+1+l+n} \\ c_{2k+l+n} \\ \vdots \\ c_{k+1+l+n} \end{pmatrix} \qquad (2.6) \\
&+ \bar{R}_{2k+1+l}(c_{k+l+n}, \dots, c_0)
\end{aligned}
$$

where we used $T^{2i}_{y_{j,i}} = 0$ if $j < 0$ for convenience of notation, i.e. in case of $k > n$ the matrix $S_k(c_{k+n}, \dots, c_0) \in \mathbb{K}^{k+1,k+1}$ is an upper band matrix with bandwith $n$. We adopt this notation from [VFZ]. Note also that $S_k(c_{k+n}, \dots, c_0)$ arises by bordering of $S_{k-1}(c_{k-1+n}, \dots, c_0) \in \mathbb{K}^{k,k}$ with column and row $k+1$.

Thus, for $m \geq 1$ fixed, the system of undetermined coefficients reads by (2.6) for $k = 0, \dots, m-1$ and $l = 0, 1$

$\boxed{k = 0,\ l = 0}$ $$T^1(c_{1+n}, \dots, c_0) = \Gamma_0(0) \cdot S_0(c_n, \dots, c_0) \cdot c_{1+n} + \bar{R}_1(c_n, \dots, c_0)$$

$\boxed{k = 0,\ l = 1}$ $$T^2(c_{2+n}, \dots, c_0) = \Gamma_0(1) \cdot S_0(c_n, \dots, c_0) \cdot c_{2+n} + \bar{R}_2(c_{1+n}, \dots, c_0)$$

$\boxed{k = 1,\ l = 0}$ $$T^3(c_{3+n}, \dots, c_0) = \Gamma_1(0) \cdot S_1(c_{1+n}, \dots, c_0) \cdot \begin{pmatrix} c_{3+n} \\ c_{2+n} \end{pmatrix} + \bar{R}_3(c_{1+n}, \dots, c_0)$$

$\boxed{k = 1,\ l = 1}$ $$T^4(c_{4+n}, \dots, c_0) = \Gamma_1(1) \cdot S_1(c_{1+n}, \dots, c_0) \cdot \begin{pmatrix} c_{4+n} \\ c_{3+n} \end{pmatrix} + \bar{R}_4(c_{2+n}, \dots, c_0) \qquad (2.7)$$

$\vdots \qquad \vdots \qquad \vdots \qquad \vdots$

$\boxed{k = m-1,\ l = 0}$ $$T^{2m-1}(c_{2m-1+n}, \dots, c_0) = \Gamma_{m-1}(0) \cdot S_{m-1}(c_{m-1+n}, \dots, c_0) \cdot \begin{pmatrix} c_{2m-1+n} \\ \vdots \\ c_{m+n} \end{pmatrix} + \bar{R}_{2m-1}(c_{m-1+n}, \dots, c_0)$$

$\boxed{k = m-1,\ l = 1}$ $$T^{2m}(c_{2m+n}, \dots, c_0) = \Gamma_{m-1}(1) \cdot S_{m-1}(c_{m-1+n}, \dots, c_0) \cdot \begin{pmatrix} c_{2m+n} \\ \vdots \\ c_{m+1+n} \end{pmatrix} + \bar{R}_{2m}(c_{m+n}, \dots, c_0).$$

Let us now assume the existence of $(c_{m+n}, \dots, c_0) \in \mathbb{K}^{m+n+1}$ with

$$S_{m-1}(c_{m-1+n}, \dots, c_0) = 0 \in \mathbb{K}^{m,m} \tag{2.8}$$

$$\bar{R}_1(c_n, \dots, c_0) = \dots = \bar{R}_{2m}(c_{m+n}, \dots, c_0) = 0. \tag{2.9}$$

Then (2.7) implies $T^1(c_{1+n}, \dots, c_0) = \dots = T^{2m}(c_{2m+n}, \dots, c_0) = 0$ with arbitrary coefficients $(c_{2m+n}, \dots, c_{m+n+1}) \in \mathbb{K}^m$ and the remaining equations of the system of undetermined coefficients read by (1.5) and (2.6) for $l \geq 0$

$$T^{2m+1+l}(c_{2m+1+l+n}, \dots, c_0) = \Gamma_m(l) \cdot S_m(c_{m+n}, \dots, c_0) \cdot \begin{pmatrix} c_{2m+1+l+n} \\ \vdots \\ c_{m+1+l+n} \end{pmatrix} + \bar{R}_{2m+1+l}(c_{m+l+n}, \dots, c_0)$$

$$= \begin{bmatrix} \gamma_0^{2m+l} & \cdots & \gamma_m^{2m+l} \end{bmatrix} \cdot \begin{pmatrix} 0 & 0 & \cdots & T^0_{y_{n-m,0}} \\ \vdots & \vdots & \cdots & \vdots \\ 0 & 0 & \cdots & T^{2m}_{y_{n,m}} \end{pmatrix} \cdot \begin{pmatrix} c_{2m+1+l+n} \\ \vdots \\ c_{m+1+l+n} \end{pmatrix} + \bar{R}_{2m+1+l}(c_{m+l+n}, \dots, c_0)$$

$$= \left[ \gamma_0^{2m+l} \cdot T^0_{y_{n-m,0}} + \dots + \gamma_m^{2m+l} \cdot T^{2m}_{y_{n,m}} \right] \cdot c_{m+1+l+n} + \bar{R}_{2m+1+l}(c_{m+l+n}, \dots, c_0)$$

$$= \left[ \binom{2m+1+l}{0} \cdot T^0_{y_{n-m,0}} + \dots + \binom{2m+1+l}{m} \cdot \binom{2m}{m}^{-1} \cdot T^{2m}_{y_{n,m}} \right] \cdot c_{m+1+l+n}$$

$$+ \bar{R}_{2m+1+l}(c_{m+l+n}, \dots, c_0)$$

$$=: g(l, c_{m+n}, \dots, c_0) \cdot c_{m+1+l+n} + \bar{R}_{2m+1+l}(c_{m+l+n}, \dots, c_0) = 0. \tag{2.10}$$

Note that the coefficient $g(l, c_{m+n}, \dots, c_0)$ is a polynomial in $l$ of degree at most $m \geq 1$. Further, assume the last column of $S_m(c_{m+n}, \dots, c_0)$ to be different from zero, i.e. $\left(T^0_{y_{n-m,0}}, \dots, T^{2m}_{y_{n,m}}\right) \neq 0$.

Now, starting from $l = 0$, equation (2.10) can recursively be solved by

$$c_{m+1+l+n} = -\frac{\bar{R}_{2m+1+l}(c_{m+l+n}, \dots, c_0)}{g(l, c_{m+n}, \dots, c_0)}, \tag{2.11}$$

as long as $g(l, c_{m+n}, \dots, c_0)$ is different from zero. If $g(\bar{l}_1, c_{m+n}, \dots, c_0) = 0$ for some $\bar{l}_1 \geq 0$, i.e. $\bar{l}_1$ is an integer root of $g(l, c_{m+n}, \dots, c_0) = 0$, then by (2.10), we additionally have to require

$$\bar{R}_{2m+1+\bar{l}_1}\left(c_{m+\bar{l}_1+n}, \dots, c_0\right) = 0 \tag{2.12}$$

and $c_{m+1+\bar{l}_1+n} \in \mathbb{K}$ turns out to be a free coefficient. If a second integer root $\bar{l}_2 > \bar{l}_1$ of $g(l, c_{m+n}, \dots, c_0) = 0$ occurs, we have to assume

$$\bar{R}_{2m+1+\bar{l}_2}\left(c_{m+\bar{l}_2+n}, \dots, c_{m+1+\bar{l}_1+n}, \dots, c_0\right) = 0 \tag{2.13}$$

for at least one value of $c_{m+1+\bar{l}_1+n}$ and $c_{m+1+\bar{l}_2+n}$ becomes the free coefficient. Continuing this process up to the largest integer root $\bar{l}_p$, the following results are shown.

**Theorem 1 :** Given $(c_{m+n},\dots,c_0)\in\mathbb{K}^{m+n+1},\ m\ge 1,\ n\ge 0$ with

(i) $S_{m-1}(c_{m-1+n},\dots,c_0)=0\in\mathbb{K}^{m,m}$ *and* $S_m(c_{m+n},\dots,c_0)\neq 0\in\mathbb{K}^{m+1,m+1}$

(ii) $\bar{R}_1(c_n,\dots,c_0)=\cdots=\bar{R}_{2m}(c_{m+n},\dots,c_0)=0$

We obtain

1) If $g(l,c_{m+n},\dots,c_0)\neq 0$ for $l\ge 0,\ l\in\mathbb{N}$, then $(c_{m+n},\dots,c_0)$ can uniquely be continued by (2.11) to a power series solution of $G[\,x,y,y',\dots,y^{(n)}\,]=0$ according to

$$y(x)=c_0+c_1\cdot x+\cdots+\frac{1}{(m+n)!}c_{m+n}\cdot x^{m+n}+\sum_{i=m+1+n}^{\infty}\frac{1}{i!}c_i\cdot x^i.$$

2) If $p\ge 1$ integer solutions $0\le\bar{l}_1<\cdots<\bar{l}_p,\ p\le m$ of $g(l,c_{m+n},\dots,c_0)=0$ exist, then we require additionally (2.12) and (2.13) with $\bar{l}_1,\dots,\bar{l}_p$, implying the existence of a family of power series solutions of dimension $q$ with $1\le q\le p$.

Essentially, Theorem 1 corresponds to Theorem 4.5 in [VFZ] that represents by itself a generalization of Lemma 2.3 in [DL]. Hence, nothing new is shown by Theorem 1 (that we do not restrict to differential polynomials or to algebraically closed fields $\mathbb{K}$ is unessential).

On the other hand, the derivation of Theorem 1 is based on the general Banach space expansion formula (1.4) of Lemma 1 that may now be used to extend Theorem 1 to systems of DAE's. In this sense, we try to follow [DL], where the scalar theory is used to obtain far reaching results concerning systems of differential polynomials of ordinary and partial type.

Prior to that, we add three remarks. Examples, and in particular further results, can be found in [VFZ] and [DL].

**Remarks: 1)** The first remark serves to clarify the correspondence between Theorem 1 and Theorem 4.5 in [VFZ]. In [VFZ] differential polynomials of the form

$$F[\,x,y,y',\dots,y^{(n)}\,]=0\,,\quad n\ge 0 \tag{2.14}$$

over an algebraically closed field $\mathbb{K}$ of characteristic zero are investigated with respect to power series solutions. The ansatz $y=\sum_{i=0}^{\infty}\frac{1}{i!}c_i\,x^i,\ x\in\mathbb{K}$ is plugged into (2.14), yielding an expansion of the form $\sum_{i=0}^{\infty}\frac{1}{i!}F^{(i)}(c_{i+n},\dots,c_0)\cdot x^i=0$, implying for $i\ge 0,\ n\ge 0$ the identities

$$T^i(c_{i+n},\dots,c_0)\ =\ F^{(i)}(c_{i+n},\dots,c_0).$$

In particular,

$$T^{2m+1+l}(c_{2m+1+l+n},\dots,c_0)\ =\ F^{(2m+1+l)}(c_{2m+1+l+n},\dots,c_0)$$

and using (2.6) as well as Theorem 3.6 in [VFZ], we see

$$\Gamma_m(l)\cdot S_m(c_{m+n},\dots,c_0)\cdot\begin{pmatrix}c_{2m+1+l+n}\\ \vdots\\ c_{m+1+l+n}\end{pmatrix}+\bar{R}_{2m+1+l}(c_{m+l+n},\dots,c_0)$$

$$=\mathcal{B}_m(l)\cdot\underbrace{\mathcal{S}_m(F)(c_{m+n},\dots,c_0)}_{m\text{-th separant matrix}}\cdot\begin{pmatrix}c_{2m+1+l+n}\\ \vdots\\ c_{m+1+l+n}\end{pmatrix}+r_{n+m+l}(c_{m+l+n},\dots,c_0),$$

where $\mathcal{S}_m(F)(c_{m+n}, \dots, c_0) \in \mathbb{K}^{m+1,m+1}$ is introduced in [VFZ], denoting the $m$-th separant matrix of $F$ evaluated at $(c_{m+n}, \dots, c_0)$. Further, by choosing the green brackets equal to zero, first the identity

$$\bar{R}_{2m+1+l}(c_{m+l+n}, \dots, c_0) = r_{n+m+l}(c_{m+l+n}, \dots, c_0)$$

follows, and secondly by elementary calculations using (1.5), the definition of $\mathcal{B}_m(l)$ in [VFZ] and a uniqueness argument, we obtain

$$\Gamma_m(l) = \mathcal{B}_m(l) \cdot Diag\left[\binom{2\cdot 0}{0} \ \cdots \ \binom{2\cdot m}{m}\right]^{-1}$$

as well as the final correspondence

$$S_m(c_{m+n}, \dots, c_0) = Diag\left[\binom{2\cdot 0}{0} \ \cdots \ \binom{2\cdot m}{m}\right] \cdot \overbrace{\mathcal{S}_m(F)(c_{m+n}, \dots, c_0)}^{separant\ matrix}.$$

**2)** If the differential algebraic equation $G[x, y, y', \dots, y^{(n)}] = 0$ is perturbed according to

$$G[\,x, y, y', \dots, y^{(n)}\,] + H[\,x, y, y', \dots, y^{(n)}\,] = 0$$

with

$$H[\,x, y, y', \dots, y^{(n)}\,] = O(\,|x, y, y', \dots, y^{(n)}|^{2m+1}\,),$$

then Theorem 1 remains valid, possibly varying higher order coefficients $c_{m+n+1}, c_{m+n+2}, \dots$ by recursion (2.11). This stability result is a consequence of the fact that the assumptions of Theorem 1 only depend on the first $2m$ derivatives of $G$ at the base solution $z_0$, i.e. from

$$G^i[\,z_0\,] = G^i[\,0, c_0, \dots, c_n\,] \in L\left[\,\mathbb{K}^{(2+n)\cdot i}, \mathbb{K}\,\right] \quad with \quad i = 1, \dots, 2m\,.$$

Compare also (1.6) and see [S2], [S3] for details.

**3)** In case of $n = 0$, i.e. in the algebraic case

$$G[\,x, y\,] = 0, \tag{2.15}$$

$G \in C^\infty(\mathbb{K}^2, \mathbb{K})$, the matrix $S_k(c_k, \dots, c_0)$ in (2.6) simplifies to a diagonal matrix and assuming (i), (ii), equation (2.10) reads

$$T^{2m+1+l}(c_{m+1+l}, \dots, c_0) \tag{2.16}$$

$$= \underbrace{\overbrace{\binom{2m+1+l}{m} \cdot \binom{2m}{m}^{-1}}^{\neq 0} \cdot \overbrace{T^{2m}_{y_{0,m}}}^{\neq 0}}_{= g(l, c_m, \dots, c_0) \neq 0} \cdot c_{m+1+l} + \bar{R}_{2m+1+l}(c_{m+l}, \dots, c_0) = 0,$$

where $T^{2m}_{y_{0,m}} \neq 0$ follows from $S_m(c_m, \dots, c_0) \neq 0$. Hence, a power series solution of (2.15) exists recursively by (2.11). In particular, case 2) in Theorem 1 is not possible.

Certainly, in the algebraic case, we know more. From (i) and (ii), we deduce for arbitrary $c_{m+1} \in \mathbb{K}$ by Taylor expansion

$$G\left[\,x,\ c_0 + c_1 x + \cdots + \frac{1}{m!} c_m x^m + \frac{1}{(m+1)!} c_{m+1} x^{m+1}\,\right] \tag{2.17}$$

$$= \frac{1}{(2m+1)!}\, x^{2m+1} \cdot T^{2m+1}(c_{m+1}, \dots, c_0) \;+\; x^{2m+2} \cdot \mathcal{R}_{2m+1}(x, c_{m+1}, \dots, c_0)$$

with smooth remainder function $\mathcal{R}_{2m+1}(\cdot)$. Splitting off $x^{2m+1}$ and using (2.16) with $l = 0$, we obtain the equation

$$\frac{1}{(2m+1)!} \cdot \Bigg[ \overbrace{\underbrace{\binom{2m+1}{m} \cdot \binom{2m}{m}^{-1} \cdot T^{2m}_{y_{0,m}}}_{= g(0, c_m, \dots, c_0) \neq 0} \cdot c_{m+1} + \bar{R}_{2m+1}(c_m, \dots, c_0)}^{= T^{2m+1}(c_{m+1}, \dots, c_0)} \Bigg]$$

$$+ \; x \cdot \mathcal{R}_{2m+1}(x, c_{m+1}, \dots, c_0) = 0 \tag{2.18}$$

which can obviously be solved with respect to $c_{m+1}$ by use of the implicit function theorem, i.e. a smooth function $c_{m+1}(x)$ exists with

$$c_{m+1}(0) = -\frac{\bar{R}_{2m+1}(c_{m+n}, \dots, c_0)}{g(0, c_{m+n}, \dots, c_0)}$$

and

$$G\left[ x,\, c_0 + c_1 x + \dots + \frac{1}{m!} c_m x^m + \frac{1}{(m+1)!}\, c_{m+1}(x)\, x^{m+1} \right] = 0. \tag{2.19}$$

This result merely represents a simple version of Tougeron's implicit function theorem. More precisely, by (2.17) we see

$$G\Bigg[ x,\, \underbrace{c_0 + c_1 x + \dots + \frac{1}{m!} c_m x^m}_{=: \bar{y}(x)} + \frac{1}{(m+1)!} c_{m+1} x^{m+1} \Bigg] = O(\, |x|^{2m+1} \,) \tag{2.20}$$

and by chain rule and (2.18), the identities

$$\frac{d}{dc_{m+1}} G\left[ x, \bar{y}(x) + \frac{1}{(m+1)!} c_{m+1} x^{m+1} \right]$$

$$= G_y\left[ x, \bar{y}(x) + \frac{1}{(m+1)!} c_{m+1} x^{m+1} \right] \cdot \frac{1}{(m+1)!}\, x^{m+1}$$

$$\overset{(2.17)}{\underset{(2.18)}{=}} \frac{1}{(2m+1)!}\, x^{2m+1} \cdot \underbrace{g(0, c_m, \dots, c_0)}_{\neq 0} \;+\; x^{2m+2} \cdot \frac{d}{dc_{m+1}} \mathcal{R}_{2m+1}(x, c_{m+1}, \dots, c_0)$$

are valid, implying under consideration of (2.20) and $c_{m+1} = 0$

$$G_y[\, x, \bar{y}(x) \,] = O(\, |x|^m \,) = \; x^m \cdot \overbrace{s(x, c_{m+1}, \dots, c_0)}^{\neq 0}$$

$$G[\, x, \bar{y}(x) \,] = O(\, |x|^{2m+1} \,)$$

with smooth remainder function $s(\cdot) \neq 0$. Thus, $\bar{y}(x)$ is an approximative solution curve of $G[x, y] = 0$ of order $2m + 1$ with $y$-derivative $G_y[x, \bar{y}(x)]$ varying at most by order $m$, i.e. the $y$-derivative $G_y$ changes quickly along $\bar{y}(x)$ near $x = 0$ compared to the change of $G$ along $\bar{y}(x)$. Then, by Tougeron's implicit function theorem [H], [R], a smooth solution curve $y(x)$ exists, agreeing with $\bar{y}(x)$ by order of $m + 1$, i.e. we end up with $G[x, y(x)] = 0$ and

$$y(x) - \bar{y}(x) = O(\, |x|^{\,m+1} \,),$$

repeating (2.19), as claimed above.

Summarizing, in the algebraic limit $n=0$, Theorem 1 turns into a 1-dimensional version of Tougeron's implicit function theorem. Now both, Tougeron's implicit function theorem as well as Lemma 1 remain valid under rather general constellations, suggesting to extend Theorem 1 to systems of DAE's.

### 3. Nonlinear systems of DAE's

In this section, Lemma 1 is used to derive sufficient conditions for the existence of power series solutions of systems of ordinary DAE's of order $n \geq 0$ given by

$$G\left[\, x, y, y', \dots, y^{(n)} \,\right] = 0 \tag{3.1}$$

with $x \in \mathbb{K}$, $y, y', \dots, y^{(n)} \in \mathbb{K}^d$ and $G \in C^\infty\left(\mathbb{K}^{1+(1+n)\cdot d}, \mathbb{K}^r\right)$, $r \geq 1$, i.e. we aim to construct a power series vector $y(x) \in \mathbb{K}^d$ satisfying $r \geq 1$ differential algebraic equations.

In principle, formulas (2.4)-(2.10) only have to be adapted to systems. We repeat (2.4) for systems according to

$$z(\varepsilon) = \sum_{i=0}^{\infty} \frac{1}{i!}\, \varepsilon^i \cdot \underbrace{\begin{pmatrix} x_i \\ y_{0,i} \\ \vdots \\ y_{n,i} \end{pmatrix}}_{=z_i} = \underbrace{\begin{pmatrix} 0 \\ c_0 \\ \vdots \\ c_n \end{pmatrix}}_{=z_0} + \varepsilon \cdot \underbrace{\begin{pmatrix} 1 \\ c_1 \\ \vdots \\ c_{n+1} \end{pmatrix}}_{=z_1} + \frac{1}{2}\varepsilon^2 \cdot \underbrace{\begin{pmatrix} 0 \\ c_2 \\ \vdots \\ c_{n+2} \end{pmatrix}}_{=z_2} + \cdots,$$

now with $y_{0,i}, \dots, y_{n,i} \in \mathbb{K}^d$, $z_i \in \mathbb{K}^{1+(1+n)\cdot d}$, $c_i \in \mathbb{K}^d$ and corresponding initial value of (3.1) by

$$G[\, z_0 \,] = G[\, 0, c_0, \dots, c_n \,] = 0 \in \mathbb{K}^r.$$

Then, by Lemma 1 and (2.2), we obtain for $k \geq 0$, $l \geq 0$,

$$T^{2k+1+l}(c_{2k+1+l+n}, \dots, c_0)$$

$$= \underbrace{\left[\boxed{T^0_{x_0}} \cdots \boxed{T^{2k}_{x_k}}\right] \cdot \boxed{C^{2k+l}} \cdot \begin{pmatrix} x_{2k+1+l} \\ x_{2k+l} \\ \vdots \\ x_{k+1+l} \end{pmatrix}}_{= \begin{cases} T^0_{x_0} & if \;\; k=l=0 \\ 0 \in \mathbb{K}^r & otherwise \end{cases}} + \left[\boxed{T^0_{y_{0,0}}} \cdots \boxed{T^{2k}_{y_{0,k}}}\right] \cdot \boxed{C^{2k+l}} \cdot \begin{pmatrix} c_{2k+1+l} \\ c_{2k+l} \\ \vdots \\ c_{k+1+l} \end{pmatrix} \tag{3.2}$$

$$+ \cdots + \left[\boxed{T^0_{y_{n,0}}} \cdots \boxed{T^{2k}_{y_{n,k}}}\right] \cdot \boxed{C^{2k+l}} \cdot \begin{pmatrix} c_{2k+1+l+n} \\ c_{2k+l+n} \\ \vdots \\ c_{k+1+l+n} \end{pmatrix} + R_{2k+1+l}(c_{k+l+n}, \dots, c_0)$$

with $T^{2i}_{x_i} \in \mathbb{K}^r$, $T^{2i}_{y_{j,i}} \in \mathbb{K}^{r,d}$, $i = 0, \dots, k$, $j = 0, \dots, n$ and $R_{2k+1+l}(c_{k+l+n}, \dots, c_0) \in \mathbb{K}^r$. As before, combining same coefficient vectors implies

$$T^{2k+1+l}(c_{2k+1+l+n},\dots,c_0)$$

$$= \underbrace{\left[\gamma_0^{2k+l} \;\cdots\; \gamma_k^{2k+l}\right]}_{=\Gamma_k(l)} \cdot \underbrace{\begin{pmatrix} T^0_{y_{n,0}} & T^0_{y_{n-1,0}} & \cdots & T^0_{y_{n-k,0}} \\ & T^2_{y_{n,1}} & \cdots & T^2_{y_{n-k+1,1}} \\ & & \cdots & T^4_{y_{n-k+2,2}} \\ & & \ddots & \vdots \\ & & & T^{2k}_{y_{n,k}} \end{pmatrix}}_{=S_k(c_{k+n},\dots,c_0)} \cdot \begin{pmatrix} c_{2k+1+l+n} \\ c_{2k+l+n} \\ \vdots \\ c_{k+1+l+n} \end{pmatrix} \tag{3.3}$$

$$+\ \bar{R}_{2k+1+l}(c_{k+l+n},\dots,c_0)$$

with $S_k(c_{k+n},\dots,c_0) \in \mathbb{K}^{(1+k)\cdot r,(1+k)\cdot d}$ and $\gamma_i^{2i+l}$, $i = 0,\dots,k$ multiplying vectors of dimension $r$. Again, we used $T^{2i}_{y_{j,i}} = 0 \in \mathbb{K}^{r,d}$ if $j < 0$, and as before note that $S_k(c_{k+n},\dots,c_0) \in \mathbb{K}^{(1+k)\cdot r,(1+k)\cdot d}$ arises by bordering of $S_{k-1}(c_{k-1+n},\dots,c_0) \in \mathbb{K}^{k\cdot r,k\cdot d}$ appropriately.

Without change of notation, formulas (2.7) remain valid for $m \geq 1$ fixed, and analogously to (2.8), (2.9), we assume the existence of $(c_{m+n},\dots,c_0) \in \mathbb{K}^{(m+n+1)\cdot d}$ with

$$S_{m-1}(c_{m-1+n},\dots,c_0) = 0 \in \mathbb{K}^{m\cdot r,m\cdot d}$$

$$\bar{R}_1(c_n,\dots,c_0) = \cdots = \bar{R}_{2m}(c_{m+n},\dots,c_0) = 0 \in \mathbb{K}^r,$$

implying $T^1(c_{1+n},\dots,c_0) = \cdots = T^{2m}(c_{2m+n},\dots,c_0) = 0 \in \mathbb{K}^r$ with arbitrary coefficients $(c_{2m+n},\dots,c_{m+n+1}) \in \mathbb{K}^{m\cdot d}$. Further, the remaining equations of the system of undetermined coefficients read by (1.5) and (3.3) for $l \geq 0$

$$T^{2m+1+l}(c_{2m+1+l+n},\dots,c_0) = \Gamma_m(l)\cdot S_m(c_{m+n},\dots,c_0)\cdot \begin{pmatrix} c_{2m+1+l+n} \\ \vdots \\ c_{m+1+l+n} \end{pmatrix} + \bar{R}_{2m+1+l}(c_{m+l+n},\dots,c_0)$$

$$= \left[\gamma_0^{2m+l} \;\cdots\; \gamma_m^{2m+l}\right] \cdot \begin{pmatrix} 0 & 0 & \cdots & T^0_{y_{n-m,0}} \\ \vdots & \vdots & \cdots & \vdots \\ 0 & 0 & \cdots & T^{2m}_{y_{n,m}} \end{pmatrix} \cdot \begin{pmatrix} c_{2m+1+l+n} \\ \vdots \\ c_{m+1+l+n} \end{pmatrix} + \bar{R}_{2m+1+l}(c_{m+l+n},\dots,c_0)$$

$$= \left[\gamma_0^{2m+l}\cdot \boxed{T^0_{y_{n-m,0}}} + \cdots + \gamma_m^{2m+l}\cdot \boxed{T^{2m}_{y_{n,m}}}\right]\cdot c_{m+1+l+n} + \bar{R}_{2m+1+l}(c_{m+l+n},\dots,c_0)$$

$$=: H(l,c_{m+n},\dots,c_0)\cdot c_{m+1+l+n} + \bar{R}_{2m+1+l}(c_{m+l+n},\dots,c_0) = 0 \tag{3.4}$$

with matrix $H(l,c_{m+n},\dots,c_0) \in \mathbb{K}^{r,d}$ and $\bar{R}_{2m+1+l}(c_{m+l+n},\dots,c_0) \in \mathbb{K}^d$. Our aim is to solve (3.4) recursively with respect to $c_{m+1+l+n}$ and $l \geq 0$.

In the quadratic case $r = d$, setting

$$g(l,c_{m+n},\dots,c_0) := det[\,H(l,c_{m+n},\dots,c_0)\,]\,, \tag{3.5}$$

a polynomial $g(l,c_{m+n},\dots,c_0)$ with respect to $l$ of degree at most $md \geq 1$ occurs, implying

$$c_{m+1+l+n} = -H(l, c_{m+n}, \dots, c_0)^{-1} \cdot \bar{R}_{2m+1+l}(c_{m+l+n}, \dots, c_0) \tag{3.6}$$

for $l \geq 0$, as long as $g(l, c_{m+n}, \dots, c_0)$ is different from zero. If $g(\bar{l}_1, c_{m+n}, \dots, c_0) = 0$ for some $\bar{l}_1 \geq 0$, then by (3.4) we additionally have to require

$$\bar{R}_{2m+1+\bar{l}_1}(c_{m+\bar{l}_1+n}, \dots, c_0) \in R[\, H(\bar{l}_1, c_{m+n}, \dots, c_0)\,] \tag{3.7}$$

with $R[\,\cdot\,]$ denoting the range of the matrix in brackets. An affine subspace of solutions arises with dimension given by the nullspace of $H(\bar{l}_1, c_{m+n}, \dots, c_0)$. If a second integer root $\bar{l}_2 > \bar{l}_1$ of $g(l, c_{m+n}, \dots, c_0) = 0$ occurs, we have to assume

$$\bar{R}_{2m+1+\bar{l}_2}(c_{m+\bar{l}_2+n}, \dots, c_{m+1+\bar{l}_1+n}, \dots, c_0) \in R[\, H(\bar{l}_2, c_{m+n}, \dots, c_0)\,] \tag{3.8}$$

for at least one value of $c_{m+1+\bar{l}_1+n}$ and a second affine subspace arises. Continuing this process up to the largest integer root $\bar{l}_p$ of $g(l, c_{m+n}, \dots, c_0) = 0$, the following results are shown in the quadratic case $r = d$.

**Theorem 2 :** Given $(c_{m+n}, \dots, c_0) \in \mathbb{K}^{(m+n+1)\cdot d},\ m \geq 1,\ d \geq 1,\ n \geq 0$ with

(i) $S_{m-1}(c_{m-1+n}, \dots, c_0) = 0 \in \mathbb{K}^{m\cdot d, m\cdot d}$ *and* $S_m(c_{m+n}, \dots, c_0) \neq 0 \in \mathbb{K}^{(m+1)\cdot d, (m+1)\cdot d}$

(ii) $\bar{R}_1(c_n, \dots, c_0) = \cdots = \bar{R}_{2m}(c_{m+n}, \dots, c_0) = 0 \in \mathbb{K}^d$

We obtain

1) If $g(l, c_{m+n}, \dots, c_0) \neq 0$ for $l \geq 0,\ l \in \mathbb{N}$, then $(c_{m+n}, \dots, c_0)$ can uniquely be continued by (3.6) to a power series solution of $G[\, x, y, y', \dots, y^{(n)}\,] = 0$ according to

$$y(x) = c_0 + c_1 \cdot x + \cdots + \frac{1}{(m+n)!}\, c_{m+n} \cdot x^{m+n} + \sum_{i=m+1+n}^{\infty} \frac{1}{i!}\, c_i \cdot x^i.$$

2) If $p \geq 1$ integer solutions $0 \leq \bar{l}_1 < \cdots < \bar{l}_p\,,\ p \leq md$ of $g(l, c_{m+n}, \dots, c_0) = 0$ exist, then we require additionally (3.7) and (3.8) with $\bar{l}_1, \dots, \bar{l}_p$, implying the existence of a family of power series solutions of dimension $q$ with $1 \leq q \leq md^2$.

**Remarks : 1)** In the underdetermined case $r < d$, a similar theorem is obviously valid, where the determinant in (3.5) has to be replaced by the set of $(r \times r)$-minors of the matrix $H(l, c_{m+n}, \dots, c_0) \in \mathbb{K}^{r,d}$. Then, if for every $l \geq 0$, a minor different from zero exists, i.e. if the rank of the matrix equals $r$ for $l \geq 0, l \in \mathbb{N}$, we obtain a family of power series solutions. If a rank drop occurs at some $l = \bar{l}_1$, then, as usual, solvability has to be required by (3.7), (3.8) for continuing the procedure.

**2)** As in the scalar case, Theorem 2 remains valid under perturbations of $G$ by order of $2m+1$.

**3)** In the algebraic case $n = 0$ with

$$G[\, x, y\,] = 0,$$

$G \in C^\infty(\mathbb{K}^{1+d}, \mathbb{K}^d),\ d \geq 1$, the matrix $S_k(c_k, \dots, c_0)$ in (3.3) simplifies again to a diagonal matrix and assuming (i), (ii), equation (3.4) reads

$$T^{2m+1+l}(c_{m+1+l}, \dots, c_0) = \overbrace{\gamma_m^{2m+l}}^{\neq 0} \cdot \boxed{T^{2m}_{y_{0,m}}} \cdot c_{m+1+l} + \bar{R}_{2m+1+l}(c_{m+l}, \dots, c_0) = 0. \tag{3.9}$$

If $det(T^{2m}_{y_{0,m}}) \neq 0$, then equation (3.9) can uniquely be solved with respect to $c_{m+1+l}$ implying the existence of a power series solution according to case 1) of Theorem 2.

Further, from (i) and (ii), we deduce for arbitrary $c_{m+1} \in \mathbb{K}^d$ by Taylor expansion

$$G\left[x,\ c_0 + c_1 x + \cdots + \frac{1}{m!} c_m x^m + \frac{1}{(m+1)!} c_{m+1} x^{m+1}\right] \tag{3.10}$$

$$= \frac{1}{(2m+1)!}\ x^{2m+1} \cdot T^{2m+1}(c_{m+1}, \dots, c_0)\ +\ x^{2m+2} \cdot \mathcal{R}_{2m+1}(x, c_{m+1}, \dots, c_0)$$

and splitting off $x^{2m+1}$ as well as using (3.9) with $l = 0$, we obtain the equation

$$\frac{1}{(2m+1)!} \cdot \left[\gamma_m^{2m} \cdot \boxed{T^{2m}_{y_{0,m}}} \cdot c_{m+1} + \bar{R}_{2m+1}(c_m, \dots, c_0)\right]\ +\ x \cdot \mathcal{R}_{2m+1}(x, c_{m+1}, \dots, c_0) = 0 \tag{3.11}$$

which can obviously be solved with respect to $c_{m+1}$ by use of the implicit function theorem using $det(T^{2m}_{y_{0,m}}) \neq 0$, i.e. a smooth function $c_{m+1}(x)$ exists with

$$G\left[x,\ c_0 + c_1 x + \cdots + \frac{1}{m!} c_m x^m\ +\ \frac{1}{(m+1)!}\ c_{m+1}(x)\ x^{m+1}\right] = 0. \tag{3.12}$$

Now, in contrast to the scalar situation with $d = 1$, in the next step we will see that the existence of $c_{m+1}(x)$ in (3.12) cannot be ascertained by Tougeron's implicit function theorem for $d \geq 2$.

First, note that

$$G[\,x,\ \underbrace{c_0 + c_1 x + \cdots + \frac{1}{m!} c_m x^m}_{=\bar{y}(x)}\ +\ \frac{1}{(m+1)!} c_{m+1} x^{m+1}\,]\ =\ O(\,|x|^{2m+1}\,) \tag{3.13}$$

by (3.10), and secondly by chain rule and (3.11), the identities

$$\frac{d}{dc_{m+1}} G\left[x, \bar{y}(x) + \frac{1}{(m+1)!} c_{m+1} x^{m+1}\right]$$

$$= G_y\left[x, \bar{y}(x) + \frac{1}{(m+1)!} c_{m+1} x^{m+1}\right] \cdot \frac{1}{(m+1)!}\ x^{m+1}$$

$$\overset{(3.10)}{\underset{(3.11)}{=}} \frac{1}{(2m+1)!}\ x^{2m+1} \cdot \gamma_m^{2m} \cdot \boxed{T^{2m}_{y_{0,m}}}\ +\ x^{2m+2} \cdot \frac{d}{dc_{m+1}} \mathcal{R}_{2m+1}(x, c_{m+1}, \dots, c_0)$$

are valid, implying

$$G_y[\,x, \bar{y}(x)\,] =\ x^m \cdot A(x)$$

$$det\{\,G_y[\,x, \bar{y}(x)\,]\,\} =\ x^{m \cdot d} \cdot \overbrace{det\,\{\,A(x)\,\}}^{\neq 0}$$

with $c_{m+1} = 0$ and $A(x) \in \mathbb{K}^{d,d}$. Hence, the determinant of the Jacobian with respect to $y$ varies along the approximative solution curve $\bar{y}(x) = c_0 + c_1 x + \cdots + \frac{1}{m!} c_m x^m$ by order of $md$ and we need

$$G[\,x, \bar{y}(x)\,] =\ O(\,|x|^{2md}\,)$$

for application of Tougeron's implicit function theorem (cf. [H], [R] with $c = 0$). Now, by (3.13), we only have $G[x, \bar{y}(x)] = O(|x|^{2m+1})$ and due to $2md > 2m+1$ for $m \geq 1,\ d \geq 2$, the assumptions of Tougeron's implicit function theorem are not satisfied along $\bar{y}(x)$.

In other words, for application of Tougeron's implicit function theorem, the approximation $\bar{y}(x)$ has to be refined to higher order. This is well known, that measuring the degeneracy of a solution curve in terms of the determinant of the Jacobi matrix along the approximation may not be optimal. A refinement of the measure seems to be given by choosing appropriate subspaces with corresponding filtration and projection operators, as shown in [ELG] and [BK], [S3].

Anyway, even for the simple constellation of Theorem 2 with $G[x, y] = 0,\ G \in C^\infty(\mathbb{K}^{1+d}, \mathbb{K}^d)$, we obtain a generalization of Tougeron's implicit function theorem (admitted in a specific constellation) and hence, it may be worth to state the algebraic result with $n = 0$ separately.

**Corollary 1 :** Given $(c_m, \dots, c_0) \in \mathbb{K}^{(m+1)\cdot d},\ m \geq 1,\ d \geq 1$ with

(i) $\quad T^0_{y_{0,0}} = \cdots = T^{2(m-1)}_{y_{0,m-1}} = 0 \in \mathbb{K}^{d,d} \quad and \quad det(\, T^{2m}_{y_{0,m}}\,) \neq 0$

(ii) $\quad \bar{R}_1(c_0) = \cdots = \bar{R}_{2m}(c_m, \dots, c_0) = 0 \in \mathbb{K}^d,$

then a smooth function $c_{m+1}(x)$ exists with

$$c_{m+1}(0) = -\frac{1}{\gamma_m^{2m}} \cdot \boxed{T^{2m}_{y_{0,m}}}^{-1} \cdot \bar{R}_{2m+1}(c_m, \dots, c_0)$$

$$G\left[x,\ c_0 + c_1 x + \cdots + \frac{1}{m!} c_m x^m + \frac{1}{(m+1)!}\, c_{m+1}(x)\, x^{m+1}\right] = 0.$$

*Matthias Stiefenhofer*

*University of Applied Sciences*

*87435 Kempten (Germany)*

*matthias.stiefenhofer@hs-kempten.de*